\documentclass[12pt]{amsart}
\usepackage{amssymb}
\usepackage{times}
\usepackage{epsf}
\usepackage{colordvi}

\renewcommand{\leq}{\leqslant}  
\renewcommand{\geq}{\geqslant}

\newtheorem{thm}{Theorem}

\begin{document}

\title{The number of transversals to line segments in $\mathbb{R}^3$}
\thanks{Research initiated at the Second McGill-INRIA Workshop on Computational
  Geometry in Computer Graphics, February 7--14 2003 co-organized by
  H. Everett, S. Lazard, and S. Whitesides, and held at the Bellairs
  Research Institute of McGill University in Holetown, St. James,
  Barbados, West Indies. }

 \author[H.~Br{\"o}nnimann, H.~Everett, S.~Lazard, F.~Sottile,
\and S.~Whitesides]{Herv\'e Br{\"o}nnimann \and Hazel Everett \and
Sylvain Lazard \and\\  Frank Sottile \and Sue Whitesides}

\address[Herv\'e Br{\"o}nnimann]{Polytechnic University, New York,  USA}
\email[Herv\'e Br{\"o}nnimann]{hbr@poly.edu}

\address[Hazel Everett and Sylvain Lazard]{LORIA (INRIA, U. Nancy 2), Nancy, France}
\email[Hazel Everett]{everett@loria.fr}

\email[Sylvain Lazard]{lazard@loria.fr}

\address[Frank Sottile]{U. Massachusetts, Amherst, MA, USA}
\email[Frank Sottile]{sottile@math.umass.edu}

\address[Sue Whitesides]{School of Computer Science, McGill University, Montr\'eal, Qc.,
        Canada}
\email[Sue Whitesides]{sue@cs.mcgill.ca}

\thanks{Research of Sottile supported in part by NSF grant DMS-0134860.}
\thanks{Research of Whitesides supported by NSERC research grants.}


\begin{abstract}
  We completely describe the structure of the connected components of transversals to a collection
  of $n$ line segments in $\mathbb{R}^3$.  We show that $n\geq 3$ arbitrary line segments in
  $\mathbb{R}^3$ admit $0, 1, \ldots, n$ or infinitely many line transversals.  In the latter case,
  the transversals form up to $n$ connected components.
\end{abstract}

\maketitle
\section{Introduction}

A $k$-transversal to a family of convex sets in $\mathbb{R}^d$ is an affine subspace of dimension
$k$ (e.g.\ a point, line, plane, or hyperplane) that intersects every member of the family.
Goodman, Pollack, and Wenger~\cite{gtt93} and Wenger~\cite{wenger98} provide two extensive surveys of
the rich subject of geometric transversal theory. In this paper, we are interested in
$1$-transversals (also called line transversals, or simply transversals) to line segments. In
$\mathbb{R}^2$, this topic was studied in the 1980's by Edelsbrunner at al.~\cite{EMPRWW82}; here we
study the topic in $\mathbb{R}^3$.

We address the following basic question.  What is the geometry and cardinality of
the set of transversals to an arbitrary collection of $n$ line segments in $\mathbb{R}^3$?  Here a
segment may be open, semi-open, or closed, and it may degenerate to a point; segments may intersect
or even overlap.  Since a line in $\mathbb{R}^3$ has four degrees of freedom, it can intersect at
most four lines or line segments in generic  position. Conversely,   it is well-known that four lines or line segments in
generic position admit $0$ or $2$ transversals; moreover,  $4$ arbitrary lines in
$\mathbb{R}^3$ admit $0, 1, 2$ or infinitely many transversals~\cite[p.  164]{hc-gi-52}.  In
contrast, our work shows  that $4$  arbitrary line segments admit up to $4$ or
infinitely many transversals.

Our interest in line transversals to segments in $\mathbb{R}^3$ is motivated by visibility problems.
In computer graphics and robotics, scenes are often represented as unions of not necessarily
disjoint polygonal or polyhedral objects.  The objects that can be seen in a particular direction
from a moving viewpoint may change when the line of sight becomes tangent to one or more objects in
the scene. Since the line of sight then becomes a transversal to a subset of the edges of the
polygons and polyhedra representing the scene, questions about transversals to segments arise very
naturally in this context.

As an example, the visibility complex~\cite{fredo,PV96} and its visibility skeleton~\cite{DDP97}
are data structures that encode visibility information of a scene; an edge of these structures
corresponds to a set of segments lying in line transversals to some $k$ edges of the scene.
Generically in $\mathbb{R}^3$, $k$ is equal to $3$ but in degenerate configurations $k$ can be
arbitrarily large.  
Such degenerate configurations frequently arise in realistic scenes; as an example a group of chairs
may admit infinitely many lines tangent to arbitrarily many of them.
It is thus essential for computing these data structures to
characterize and compute the transversals to $k$ segments in $\mathbb{R}^3$.  Also, to bound the
size of the visibility complex one needs to bound the number of connected components of
transversals to $k$ arbitrary line segments. While the answer $O(1)$ suffices for giving asymptotic
results, the present paper establishes the actual bound.

As mentioned above, in the context of 3D visibility, lines tangent to objects are more relevant than
transversals; lines tangent to a polygon or polyhedron along an edge happen to be
transversals to this edge.  (For bounds on the space of transversals to convex polyhedra in
$\mathbb{R}^3$ see~\cite{PS}.) The literature related to lines tangent to objects falls into two
categories. The one closest to our work deals with characterizing 
the degenerate configurations of curved objects with respect to tangent lines.
MacDonald, Pach and Theobald~\cite{theo1} give a complete description of the set of lines tangent to
four unit balls in $\mathbb{R}^3$.  Megyesi, Sottile and Theobald~\cite{MST03} describe the set of
lines meeting two lines and tangent to two spheres in~$\mathbb{R}^3$, or tangent to two quadrics in
$\mathbb{P}^3$. Megyesi and Sottile~\cite{MS03} describe the set of lines meeting one line and
tangent to two or three spheres in~$\mathbb{R}^3$.  A nice survey of these results can be found
in~\cite{theobald_habilitation}.

The other category of results deals with lines tangent to $k$ among $n$ objects in $\mathbb{R}^3$.
For polyhedral objects, De Berg, Everett and Guibas~\cite{BEG98} showed a $\Omega(n^3)$ lower bound
on the number of free (i.e., non-occluded by the interior of any object) lines tangent to $4$
amongst $n$ disjoint homothetic convex polyhedra.  Br{\"o}nnimann et al.~\cite{cccg02} showed that,
under a certain general position assumption, the number of lines tangent to $4$ amongst $k$ bounded
disjoint convex polyhedra of total complexity $n$ is $O(n^2k^2)$.  For curved objects, Devillers et
al.~\cite{DDE03} and Devillers and Ramos~\cite{DR01} (see also \cite{DDE03}) presented simple
$\Omega(n^2)$ and $\Omega(n^3)$ lower bounds on the number of free maximal
segments tangent to $4$ amongst $n$ unit balls and amongst $n$ arbitrarily sized balls. 
Agarwal, Aronov and Sharir~\cite{AAS99} showed an upper bound of $O(n^{3+\epsilon})$ on the
complexity of the space of line transversals of $n$ balls by studying the lower envelope of a set of
functions; a study of the upper envelope of the same set of functions yields the same upper bound on
the number of free lines tangent to four balls~\cite{DR01}. Durand et al.~\cite{fredo} showed an
upper bound of $O(n^{8/3})$ on the expected number of possibly occluded lines tangent to $4$ among
$n$ uniformly distributed unit balls.  Under the same model, Devillers et al.~\cite{DDE03} recently
showed a bound of $\Theta(n)$ on the number maximal free line segments tangent to $4$ among $n$
balls. 

\section{Our results}

We say that 
two transversals to a
collection of line segments are in the same {\em connected component} if and only if one can be
continuously moved into the other while remaining a transversal to the collection of line segments.
Equivalently, the two points in line space (e.g., in Pl{\"u}cker space) corresponding to the two
transversals are in the same connected component of the set of points corresponding to all the
transversals to the collection of line segments.

Our  main result is the following theorem.

\begin{thm}\label{T:main}
  A collection of $n\geq 3$ arbitrary line segments in $\mathbb{R}^3$ admits $0, 1, \ldots, n$ or
  infinitely many transversals.  In the latter case, the transversals can form any number, from $1$
  up to $n$ inclusive, of connected components.
\end{thm}

More precisely we show that, when $n\geq 4$, there can be more than $2$ transversals only if the
segments are in some degenerate configuration, namely if the $n$ segments are members of one ruling
of a hyperbolic paraboloid or a hyperboloid of one sheet, or if they are concurrent, or if they all
lie in a plane with the possible exception of a group of one or more segments that all meet that
plane at the same point.  

\begin{figure*}[t]
 \[
   \epsfysize=3.5cm\epsfbox{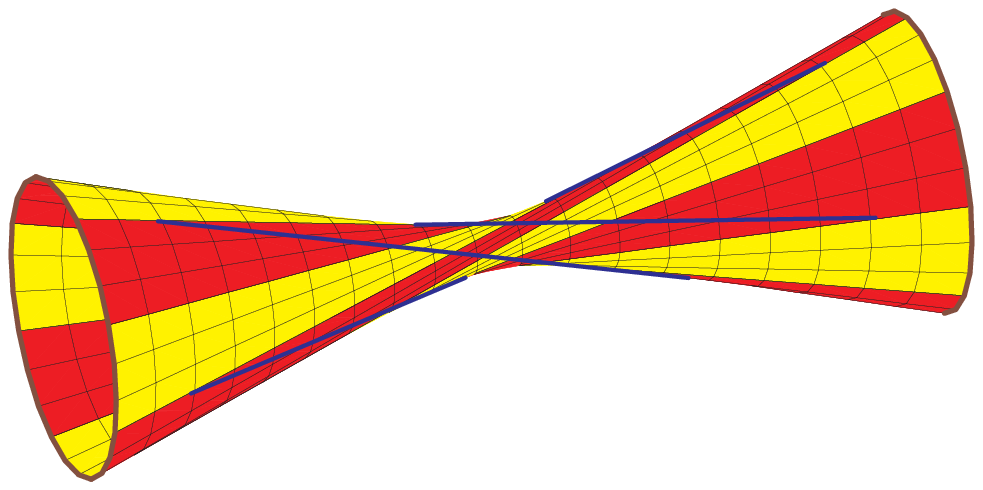}\qquad
   \epsfysize=3.5cm\epsfbox{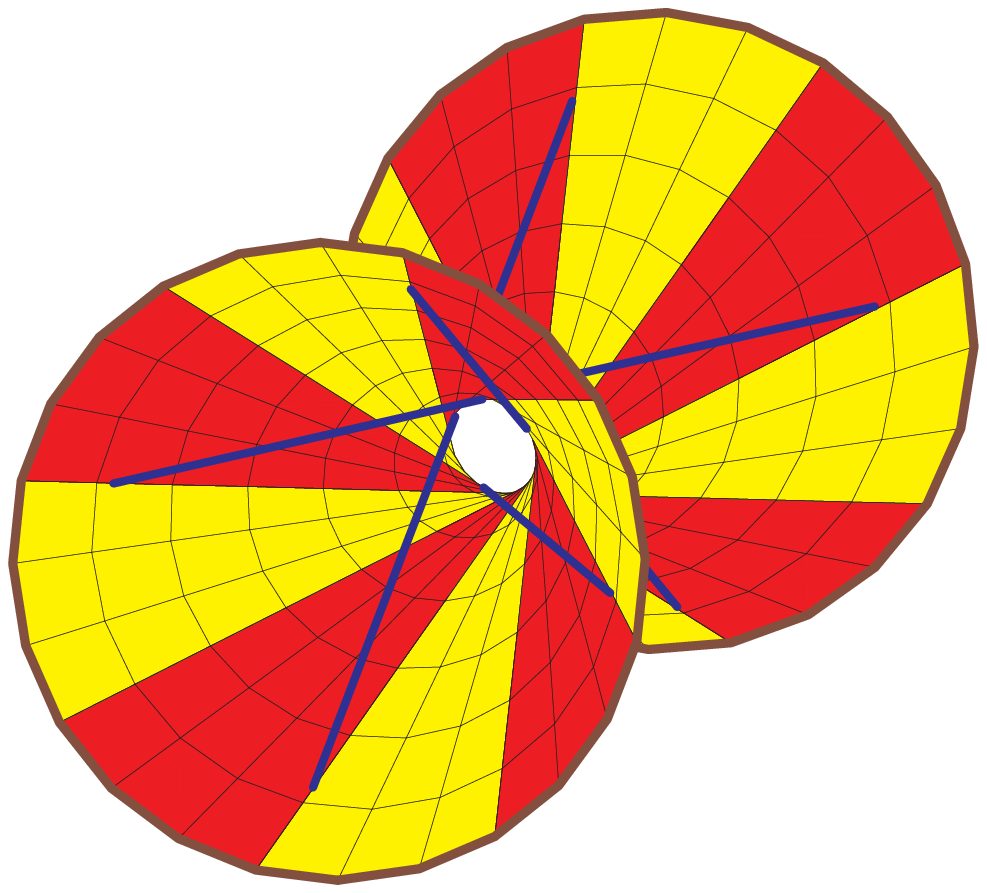}
 \]
\caption{Two views of a hyperboloid of one sheet containing four line segments and their 
  four connected components of  transversals (corresponding to the shaded regions). The four
  segments are symmetric under rotation about the axis of the hyperboloid.}
\label{F:interval}
\end{figure*}

Moreover,  in these degenerate configurations the  number of connected components of transversals is as
follows.  If the segments are members of one ruling of a hyperbolic paraboloid, or if they are
concurrent, their transversals form at most one connected component.  If they are members of one
ruling of a hyperboloid of one sheet, or if they are coplanar, their transversals can have up to $n$
connected components (see Figures~\ref{F:interval} and \ref{F:coplanar}).  Finally, if the segments
all lie in a plane with the exception of a group of one or more segments that all meet that plane at
the same point, their transversals can form up to $n-1$ connected components (see
Figures~\ref{F:triangle} and \ref{fig:example-with-3-connected-components}).

The geometry of the transversals is as follows.
We consider here $n\geq 4$ segments that are pairwise non-collinear; otherwise, as we shall see in
Section~\ref{sec:proof}, we can replace segments having the same supporting line by their common
intersection.  If the segments are members of one ruling of a hyperbolic paraboloid or a hyperboloid
of one sheet, their transversals lie in the other ruling (see Figures~\ref{F:interval} and
\ref{fig:hyperbolic-paraboloid}).  If the segments are concurrent at a point $p$, their transversals
consist of the lines through $p$ and, if the segments also lie in a plane $H$, of lines in $H$.  If
the segments consist of a group segments lying in a plane $H$ and meeting at a point $p$, together
with a group of one or more segments meeting $H$ at a point $q\neq p$ and lying in a plane $K$
containing $p$, their transversals lie in $H$ and $K$ (see
Figure~\ref{fig:example-with-3-connected-components}).  Finally, if none of the previous conditions
holds and if the segments all lie in a plane $H$ with the possible exception of a group of one or
more segments that all meet $H$ at the same point, then their transversals lie in $H$ (see
Figures~\ref{F:triangle} and \ref{F:coplanar}).

If all the $n$ segments are coplanar, the set of connected components of transversals, as well as
any one of these components, can be of linear complexity~\cite{EMPRWW82}. Otherwise we prove that
each of the connected components has constant complexity and can be represented by
an interval  on a line or on a circle  (or possibly by two intervals in the case depicted in
Figure~\ref{fig:example-with-3-connected-components}).

A connected component of transversals may be an isolated line. For example, three segments forming a
triangle and a fourth segment intersecting the interior of the triangle in one point have exactly
three transversals (Figure~\ref{F:triangle} shows a similar example with infinitely many
transversals).  Also, the four segments in Figure~\ref{F:interval} can be shortened so that the four
connected components of transversals reduce to four isolated transversals.

Finally, as discussed in the conclusion, an $O(n\log n)$-time algorithm for computing the
transversals to $n$ segments directly follows from the proof of Theorem~\ref{T:main}.

\section{Proof of Theorem~\ref{T:main}}\label{sec:proof}

Every non-degenerate line segment is contained in  its \emph{supporting line}.  
We define the supporting line of a point to be the vertical line through that point.
We prove Theorem~\ref{T:main} by considering the three following 
cases which  cover all possibilities but are not exclusive.
\begin{enumerate}
\item Three supporting lines are pairwise skew.
\item Two supporting lines are coplanar.
\item All the segments are coplanar.
\end{enumerate}

We can assume in what follows that \emph{the supporting lines are pairwise distinct.}  Indeed, if
disjoint segments have the same supporting line $\ell$, then $\ell$ is the only  transversal
to those segments, and so the set of  transversals is either empty or consists of $\ell$.  If
non-disjoint segments have the same supporting line, then any  transversal must meet the
intersection of the segments.  We can replace these overlapping segments by their common
intersection.

\subsection{Three supporting   lines are skew}
\label{sec:case1}

Three pairwise skew lines lie on a unique doubly-ruled hyperboloid, namely, a hyperbolic paraboloid
or a hyperboloid of one sheet (see the discussion in~\cite[\S 3]{PW01}). Furthermore, they are
members of one ruling, say the ``first'' ruling, and their  transversals are the lines in the
``second'' ruling that are not parallel to any of the three given skew lines.

\begin{figure}[t]
\begin{center}
  \setlength{\unitlength}{.8pt}
\hspace{.5cm}  \begin{picture}(270,160)(-5,-5)
   \put(  5,  5){\epsfysize=1.6in \epsfbox{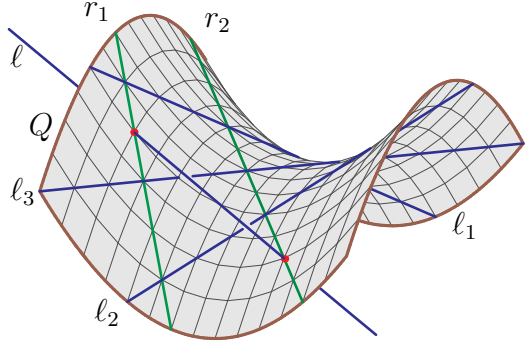}}  
   \put( 35,  8){$\ell_2$}
   \put( -5, 65){$\ell_3$}
   \put(204, 49){$\ell_1$}
   \put( -5,128){$\ell$}
   \put( 30,152){$r_1$}
   \put( 87,147){$r_2$}
   \put(  4,96){$Q$}
  \end{picture}
\caption{Line $\ell$ intersects in two points the hyperbolic paraboloid spanned by the 
  lines $\ell_1$, $\ell_2$ and $\ell_3$. The two lines $r_1$ and $r_2$ meet the four lines $\ell_1$,
  $\ell_2$, $\ell_3$, and~$\ell$.}
\label{fig:hyperbolic-paraboloid}
\end{center}
\end{figure}

Consider first the case where there exists a fourth segment whose supporting line $\ell$ does not
lie in the first ruling. Either $\ell$ is not contained in the hyperboloid or it lies in the second
ruling. In both cases, there are at most two  transversals to the four supporting lines, which
are lines of the second ruling that meet or coincide with $\ell$ (see
Figure~\ref{fig:hyperbolic-paraboloid})~\cite[p. 164]{hc-gi-52}.  Thus there are at most two 
transversals to the $n$ line segments.

Now suppose that all the $n\geq 3$ supporting lines of the segments $s_i$ lie in the first ruling of
a hyperbolic paraboloid. The lines in the second ruling can be parameterized by their intersection
points with any line $r$ of the first ruling.  Thus the set of lines in the second ruling that meet
a segment $s_i$ corresponds to an interval on line $r$.  Hence the set of transversals to the $n$
segments corresponds to the intersection of $n$ intervals on $r$, that is, to one interval on this
line, and so the  transversals form one connected component.

Consider finally the case where the $n\geq 3$ supporting lines lie in the first ruling of a
hyperboloid of one sheet (see Figure~\ref{F:interval}).  The lines in the second ruling can be
parameterized by points on a circle, for instance, by their intersection points with a circle lying
on the hyperboloid of one sheet.  Thus the set of transversals to the $n$ segments corresponds to
the intersection of $n$ intervals on this circle.  This intersection can have any number of
connected components from $0$ up to $n$ and any of these connected components may consist of an
isolated point on the circle.  The set of transversals can thus have any number of connected
components from $0$ up to $n$ and any of these connected components may consist of an isolated
transversal.  Figure~\ref{F:interval} shows two views of a configuration with $n=4$ line segments
having $4$ connected components of transversals.

In this section we have proved that if the supporting lines of $n\geq 3$ line segments lie in one
ruling of a hyperboloid of one sheet, the segments admit $0, 1, \ldots, n$ or infinitely many
transversals which form up to $n$ connected components.  If supporting lines lie in one ruling of a
hyperbolic paraboloid, the segments admit at most $1$ connected component of transversals. Otherwise
the segments admit up to $2$ transversals.

\subsection{Two supporting lines are coplanar}
\label{sec:case2}

Let $\ell_1$ and $\ell_2$ be two (distinct) coplanar supporting lines in a plane $H$.  First consider the case
where $\ell_1$ and $\ell_2$ are parallel.  Then the transversals to the $n$ segments all lie in $H$.
If some segment does not intersect $H$ then there are no transversals; otherwise, we can replace
each segment by its intersection with $H$ to obtain a set of coplanar segments, a configuration
treated in Section~\ref{sec:coplanar}.

Now suppose that $\ell_1$ and $\ell_2$ intersect at point $p$.  Consider all the supporting lines
not in $H$.  If no such line exists then all segments are coplanar; see Section~\ref{sec:coplanar}.
If such lines exist and any of them is parallel to $H$ then all transversals to the $n$ segments lie
in the plane containing $p$ and that line. We can again replace each segment by its intersection
with that plane to obtain a set of coplanar segments, a configuration treated in
Section~\ref{sec:coplanar}.

We can now assume that there exists a supporting line not in $H$.  Suppose that all the supporting
lines not in $H$ go through $p$.  If all the segments lying in these supporting lines contain $p$
then we may replace all these segments by the point $p$ without changing the set of transversals to
the $n$ segments. Then all resulting segments are coplanar, a configuration treated in
Section~\ref{sec:coplanar}. Now if some segment $s$ does not contain $p$ then the only possible
transversal to the $n$ segments is the line containing $s$ and $p$.

\begin{figure}[t]
\[
   \begin{picture}(150,110)
    \put(0,0){\epsfysize=110pt\epsfbox{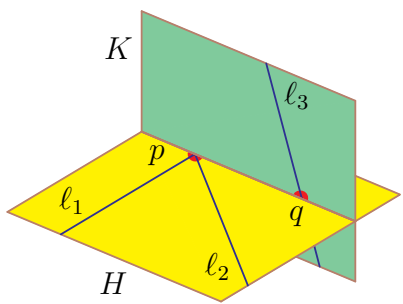}}
    \put(20,35){$\ell_1$}  \put(75,10){$\ell_2$}  \put(105,75){$\ell_3$}
    \put(54,54){$p$}  \put(107,31){$q$}
    \put(35,3){$H$}  \put(37,92){$K$}
   \end{picture}
\]
\caption{Lines $\ell_1$ and $\ell_2$ intersect at point $p$, and line
         $\ell_3$ intersects plane $H$ in a point $q$ distinct from $p$.}  
\label{fig:two-lines-meet}
\end{figure}

We can now assume that there exists a supporting line $\ell_3$ intersecting $H$ in exactly one
point $q$ distinct from $p$ (see Figure~\ref{fig:two-lines-meet}).  Let $K$ be the plane containing
$p$ and $\ell_3$.  Any  transversal to the lines $\ell_1$, $\ell_2$ and $\ell_3$ lies in $K$
and goes through $p$, or lies in $H$ and goes through $q$.

If there exists a segment $s$ that lies neither in $H$ nor in $K$ and goes through neither $p$ nor
$q$, then there are at most two transversals to the $n$ segments, namely, at most one line in $K$
through $p$ and $s$ and at most one line in $H$ through $q$ and $s$.

\begin{figure}[t]
\[
  \begin{picture}(205,120)(0,-5)
   \put(0,0){\epsfxsize=200pt\epsffile{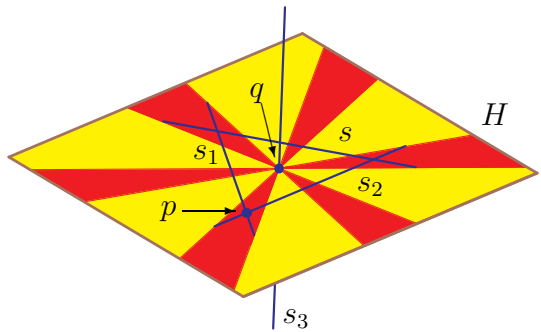}}
   \put(179,75){$H$}
   \put( 70,62){$s_1$}
   \put( 58,41){$p$} \put(66,42.3){\vector(1,0){20}}        
   \put( 91.5,85.5){$q$}\put(96,83){\vector(1,-4){5}}
   \put(132,50){$s_2$}   \put(125,69){$s$}
   \put(104, 0){$s_3$}
  \end{picture}
\]
\caption{Four segments having  three connected components of transversals.}
\label{F:triangle}
\end{figure}

We can thus assume that all segments lie in $H$ or $K$ or go through $p$ or $q$.  If there exists a
segment $s$ that goes through neither $p$ nor $q$, it lies in $H$ or $K$.  If it lies in $H$ then
all the transversals to the $n$ segments lie in $H$ (see Figure~\ref{F:triangle}).  Indeed, no line
in $K$ through $p$ intersects $s$ except possibly the line $pq$ which also lies in $H$.  We can
again replace each segment by its intersection with $H$ to obtain a set of coplanar segments; see
Section~\ref{sec:coplanar}.  The case where $s$ lies in $K$ is similar.

We can now assume that all segments go through $p$ or $q$ (or both).  Let $n_p$ be the number of
segments not containing $p$, and $n_q$ be the number of segments not containing $q$.  Note that
$n_p+n_q\leq n$.

Among the lines in $H$ through $q$, the transversals to the $n$ segments are the transversals to
the $n_{q}$ segments not containing $q$.  We can replace these $n_{q}$ segments by their
intersections with $H$ to obtain a set of $n_{q}$ coplanar segments in $H$.  The transversals to
these segments in $H$ through $q$ can form up to $n_{q}$ connected components.  Indeed, the lines in
$H$ through $q$ can be parameterized by a point on a circle, for instance, by their polar angle in
$\mathbb{R}/\pi\mathbb{Z}$. Thus the set of lines in $H$ through $q$ and through a segment in $H$
corresponds to an interval of $\mathbb{R}/\pi\mathbb{Z}$.  Hence the set of transversals to the
$n_{q}$ segments corresponds to the intersection of $n_{q}$ intervals in $\mathbb{R}/\pi\mathbb{Z}$
which can have up to $n_{q}$ connected components.

Similarly, the lines in $K$ through $p$ that are transversals to the $n$ segments can form up to
$n_{p}$ connected components.  Note furthermore that the line $pq$ is a transversal to all segments
and that the connected component of transversals that contains the line $pq$ is counted twice.
Hence there are at most $n_{p}+n_{q}-1\leq n-1$ connected components of transversals to the $n$
segments.

\begin{figure}[tb]
\[
  \begin{picture}(225,150)
   \put(4,0){\epsfxsize=225pt\epsffile{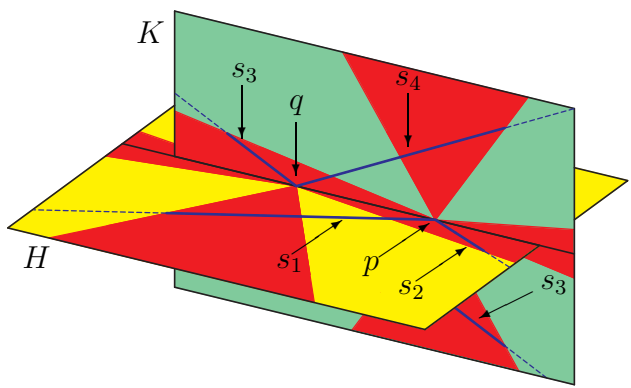}}
    \put(0,45){$H$}    \put(43,130){$K$}
    \put( 96, 45){{$s_1$}} \put(102,50){{\vector(3,2){19}}}
    \put(142, 35){{$s_2$}} \put(149,40){{\vector(3,2){19}}}
    \put(129, 44){{$p$}}   \put(135,49){{\vector(3,2){19}}}

    \put( 79,117){{$s_3$}} \put(83,114){{\vector(0,-1){20}}}
    \put(196, 37){{$s_3$}} \put(193, 36){{\vector(-2,-1){20}}}
    \put(141,114){{$s_4$}} \put(146,111){{\vector(0,-1){20}}}
    \put(101,105){{$q$}} \put(103.5,100){{\vector(0,-1){20}}}
  \end{picture}
\]
\caption{Four segments having three connected components of transversals.}
\label{fig:example-with-3-connected-components}
\end{figure}

To see that the bound of $n-1$ connected components is reached, first consider $n/2$ lines in $H$
through $p$, but not through $q$. Their transversals through $q$ are all the lines in $H$ though
$q$, except for the lines that are parallel to any of the $n/2$ given lines. This gives $n/2$
connected components. Shrinking the $n/2$ lines to sufficiently long segments still gives $n/2$
connected components of transversals in $H$ through $q$. The same construction in plane $K$ gives
$n/2$ connected components of transversals in $K$ through $p$. This gives $n-1$ connected components
of transversals to the $n$ segments since the component containing line $pq$ is counted twice.
Figure~\ref{fig:example-with-3-connected-components} shows an example of $4$ segments having $3$
connected components of transversals.

In this section we have proved that  $n\geq 3$ segments  having at least two coplanar
supporting lines either can  be reduced to $n$ coplanar segments or may
have up to $n-1$ connected components of  transversals.

\subsection{All the  line segments are  coplanar}
\label{sec:coplanar}

We prove here that $n\geq 3$ coplanar line segments in $\mathbb{R}^3$ admit up to $n$
connected components of  transversals.

Let $H$ be the plane containing all the $n$ segments.  There exists a transversal not in $H$ if and
only if all segments are concurrent at a point $p$. In this case, the transversals consist of the
lines through $p$ together with the transversals lying in $H$. To see that they form only one
connected component, notice that any transversal in $H$ can be translated to $p$ while remaining a
transversal throughout the translation.  We thus can assume in the following that all transversals
lie in $H$, and we consider the problem in $\mathbb{R}^2$.

We consider the usual geometric transform (see e.g. \cite{EMPRWW82}) where a line in $\mathbb{R}^2$
with equation $y=a x+b$ is mapped to the point $(a,b)$ in the dual space.  The transversals to a
segment are transformed to a double wedge; the double wedge degenerates to a line when the segment
is a point. The apex of the double wedge is the dual of the line containing the segment.

A transversal to the $n$ segments is represented in the dual by a point in the intersection of all
the double wedges. There are at most $n+1$ connected components of such points~\cite{EMPRWW82} (see
also~\cite[Lemma 15.3]{E87}).  Indeed, each double wedge consists of two wedges separated by the
vertical line through the apex.  The intersection of all the double wedges thus consists of at most
$n+1$ convex regions whose interiors are separated by at most $n$ vertical lines.

Notice that if there are exactly $n+1$ convex regions then two of these regions are connected at
infinity by the dual of some vertical line, in which case the segments have a vertical 
transversal. Thus the number of connected components of  transversals is at most~$n$.

\begin{figure}[t]
\[
  \begin{picture}(145,120)(0,-15)
   \put(0,0){\epsfxsize=145pt\epsffile{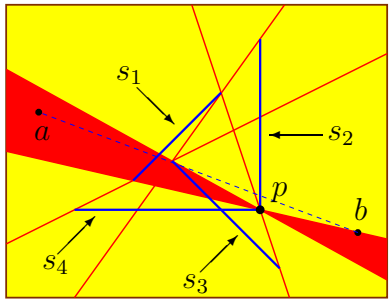}}
   \put(11,60){$a$}   \put(101,39){$p$}   \put(132,30){$b$}
    \put( 42,83){$s_1$} \put( 50,80){\vector(1,-1){14}}
    \put(122,60){$s_2$} \put(120,62){\vector(-1,0){20}}
    \put( 67, 4){$s_3$} \put(74,10){\vector(1,1){14}}
    \put( 14,11){$s_4$} \put(20,17){\vector(1,1){14}}
  \end{picture}
\]
\vspace{-10mm}
\caption{Four coplanar segments having four connected components of transversals.}
\label{F:coplanar}
\end{figure}

To see that this bound is sharp consider the configuration in Figure~\ref{F:coplanar} of $4$
segments having $4$ components of  transversals.  Three of the components consist of isolated
lines and one consists of a connected set of lines through $p$ (shaded in the figure).  Observe that
the line segment ${ab}$ meets the three isolated lines.  Thus the set of  transversals to the
four initial segments and segment $ab$ consists of the $3$ previously mentioned isolated
transversals, the line $pb$ which is isolated, and a connected set of lines through $p$.  This may
be repeated for any number of additional segments, giving configurations of $n$ coplanar line
segments with $n$ connected components of  transversals.

\section{Conclusion} 

This paper has characterized the geometry, cardinality and complexity of the set of transversals to
an arbitrary collection of line segments in $\mathbb{R}^3$.
In addition to contributing to geometric transversal theory, 
we anticipate that the results will be useful in the  
design of geometric algorithms and in their running time analyses.

While algorithmic issues have not been the main concern of the paper, we note that the proof of
Theorem~\ref{T:main} leads to an $O(n\log n)$-time algorithm in the real RAM model of computation.
First reduce in $O(n\log n$) time the set of segments to the case of pairwise distinct supporting lines.
Choose any three of these lines. Either they are pairwise skew or two of them are coplanar.  If they
are pairwise skew (see Section~\ref{sec:case1}), their transversals, and hence the transversals to
all $n$ segments, lie in one ruling of a hyperboloid. Any segment that intersects the hyperboloid in
at most two points admits at most $2$ transversals that lie in that ruling. Checking whether these
lines are transversals to the $n$ segments can be done in linear time. Consider now the case of a
segment that lies on the hyperboloid.  Its set of transversals, lying in the ruling, can be
parameterized in constant time by an interval  on a line or a circle depending on the type of the
hyperboloid. Computing the transversals to the $n$ segments thus reduces in linear time to
intersecting $n$ intervals  on a line or on a circle, which can be done in $O(n\log n)$ time.  If two
supporting lines are coplanar (see Section~\ref{sec:case2}), computing the transversals to the $n$
segments reduces in linear time to computing transversals to at most $n$  segments in one or two
planes, which can be done in $O(n\log n)$ time~\cite{EMPRWW82}.

Finally, note that we did not consider in this paper, for simplicity of the exposition, lines or
half-lines although our theorem holds when such lines in $\mathbb{R}^3$ are allowed. Note for
example that, in $\mathbb{R}^3$, the transversals to $n\geq 3$ lines of one ruling of a hyperboloid
of one sheet are all the lines of the other ruling with the exception of the lines parallel to the
$n$ given lines. Thus, in $\mathbb{R}^3$, the transversals form $n$ connected components. Remark
however that our theorem does not hold for lines in projective space $\mathbb{P}^3$; in this case,
our proof directly yields that, if a set of lines admit infinitely many transversals, they form one
connected component.

\section*{Acknowledgments}
We would like to thank the other participants of the workshop for useful discussions.

\providecommand{\bysame}{\leavevmode\hbox to3em{\hrulefill}\thinspace}
\providecommand{\MR}{\relax\ifhmode\unskip\space\fi MR }
\providecommand{\MRhref}[2]{%
  \href{http://www.ams.org/mathscinet-getitem?mr=#1}{#2}
}
\providecommand{\href}[2]{#2}

\def\baselinestretch{0.93}

\medskip



\normalsize
\small

\begin{thebibliography}{10}
\bibitem {AAS99}
P. K. Agarwal, B. Aronov, and M. Sharir.
\newblock Line transversals of balls and smallest enclosing cylinders in three dimensions.
\newblock \emph{Discrete Comput. Geom.}, 21:373--388, 1999.


\bibitem {BEG98}
M. de Berg, H. Everett, and L.J. Guibas. 
\newblock The union of moving polygonal pseudodiscs -- Combinatorial bounds and applications. 
\newblock \emph{Computational Geometry: Theory and Applications},
11:69--82, 1998. 

\bibitem{cccg02}
H.~Br{\"o}nnimann, O.~Devillers, V.~Dujmovic, H.~Everett, M.~Glisse,
X.~Goaoc, S.~Lazard, H.-S.~Na, and S.~Whitesides.
\newblock On the number of lines tangent to four convex polyhedra.
\newblock \emph{Proc. Fourteenth Canad. Conf. 
Comput. Geom.}, pp.\ 113-117, 2002.



\bibitem{DDE03}
O. Devillers, V. Dujmovi{\'c}, H. Everett, X. Goaoc, S. Lazard, H.-S. Na, and S. Petitjean.
\newblock The expected number of 3D visibility events is linear.
\newblock  To appear in \emph{SIAM J. Computing}.

\bibitem{DR01}
O. Devillers and P. Ramos.
\newblock Personal communication, 2001.

\bibitem{DDP97}
{F.~Durand, G.~Drettakis, and C.~Puech}.
\newblock {The visibility skeleton: a powerful and efficient multi-purpose  global visibility tool}.
\newblock  \emph{Proc. 31st COMPUTER GRAPHICS Ann. Conference Series 
  (SIGGRAPH'97)}, pp. 89--100, 1997.

\bibitem{fredo}
F.~Durand, G.~Drettakis, and C.~Puech.
\newblock The 3D visibility complex.
\newblock \emph{ACM Trans. Graphics} 21(2):176--206, 2002.


\bibitem{E87}
H. Edelsbrunner.
\newblock \emph{Algorithms in Combinatorial Geometry.}
\newblock   Springer-Verlag, Heidelberg, Germany, 1987. 

\bibitem{EMPRWW82}
H.~Edelsbrunner, H.~A.~Maurer, F.~P.~Preparata, A.~L.~Rosenberg, E.~Welzl, and D.~Wood.
\newblock Stabbing line segments.
\newblock \emph{BIT} 22:274-281, 1982. (See also~\cite[Chapter 15]{E87}.)


\bibitem{gtt93}
J.~E. Goodman, R.~Pollack, and R.~Wenger.
\newblock Geometric transversal theory.
\newblock In \emph{New Trends in Discrete and Computational Geometry}, (J. Pach,
ed.), Springer Verlag, Heidelberg, pp.\ 163-198, 1993.

\bibitem{hc-gi-52}
D.~Hilbert and S.~Cohn-Vossen.
\newblock \emph{Geometry and the Imagination}.
\newblock Chelsea Publishing Company, New York, 1952.


\bibitem{theo1}
I.~G. Macdonald, J.~Pach, and T.~Theobald.
\newblock Common tangents to four unit balls in $\mathbb{R}^3$.
\newblock \emph{Discrete Comput. Geom.} 26(1):1-17, 2001.
   
\bibitem{MS03}
G. Megyesi and F. Sottile.
\newblock The Envelope of Lines Meeting a fixed line and Tangent to Two Spheres. 
\newblock Manuscript. 

\bibitem{MST03}
G. Megyesi, F. Sottile, and T. Theobald.
\newblock   Common transversals and tangents to two lines and two quadrics in $\mathbb{P}^3$.
\newblock To appear in \emph{Discrete Comput. Geom.} 

      
\bibitem {PS}
{M. Pelligrini and P. W. Shor}
\newblock Finding stabbing lines in 3-space.
\newblock  \emph{Discrete Comput. Geom.}, 8:191-208, 1992.


\bibitem {PV96}
{M. Pocchiola and G. Vegter}.
\newblock   {The visibility complex}.
\newblock \emph{Internat. J. Comput. Geom. Appl.}, 6(3):279--308, 1996.

\bibitem{PW01}
H.~Pottman and J.~Wallner.
\newblock \emph{Computational Line Geometry}.
\newblock Springer-Verlag, Berlin, 2001.

\bibitem {theobald_habilitation}
T.~Theobald.
\newblock   New algebraic methods in computational geometry.
\newblock Habilitationsschrift, Fakult{\"a}t fur Mathematik, Technische Universit{\"a}t M{\"u}nchen,
Germany, 2003. 

\bibitem{wenger98}
R.~Wenger.
\newblock Progress in geometric transversal theory.
\newblock In \emph{Advances in Discrete and Computational Geometry},
  (B.~Chazelle, J.~E.~Goodman, and R.~Pollack, eds.),
   Amer. Math. Soc., Providence,  pp.\ 375--393, 1998.

\end{thebibliography}
\end{document}